\documentstyle[11pt,epsf,amstex]{article}

\def\fgCarterCurve{
\begin{figure}[htb]
\begin{center}
\mbox{\epsfxsize=2.5in \epsfbox{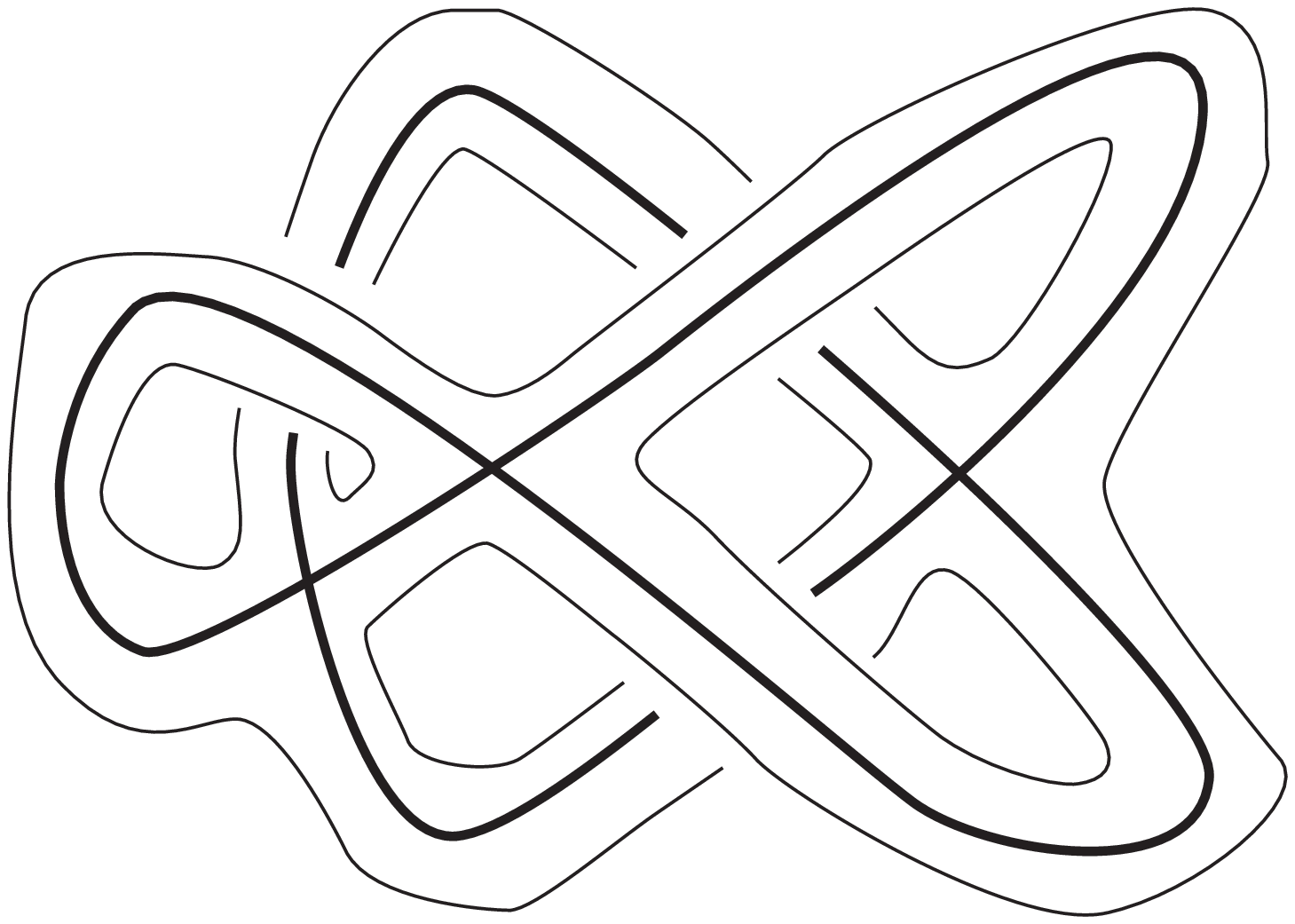}}
\end{center}
\caption{}
\label{CarterCurve}
\end{figure}}

\def\fghopfboth{
\begin{figure}[htb]
\begin{center}
\mbox{\epsfxsize=3.0in \epsfbox{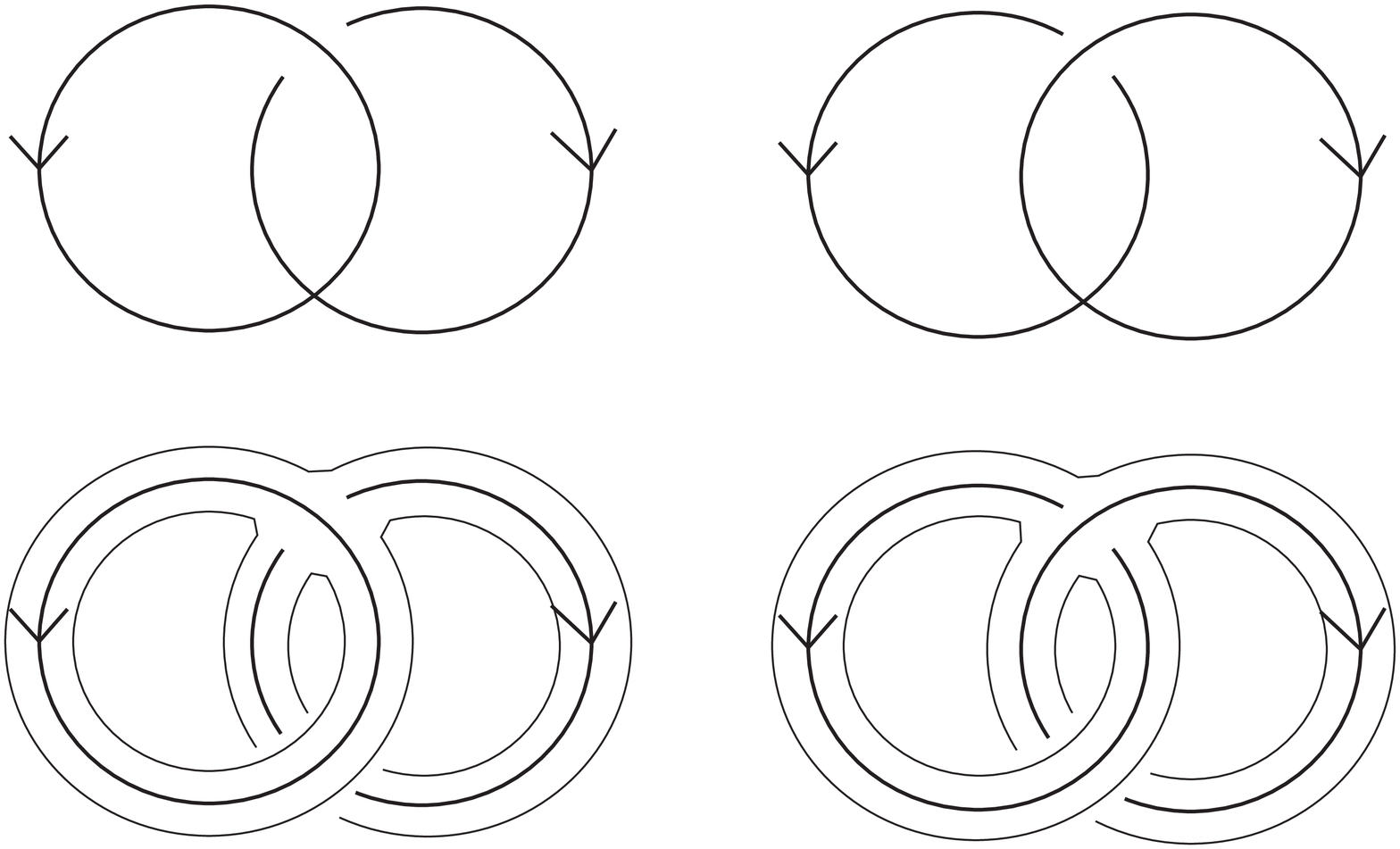}}
\end{center}
\caption{Virtual and abstract pseudo-Hopf links}
\label{hopfboth}
\end{figure}}

\def\fghopfcancel{
\begin{figure}[htb]
\begin{center}
\mbox{\epsfxsize=3.0in \epsfbox{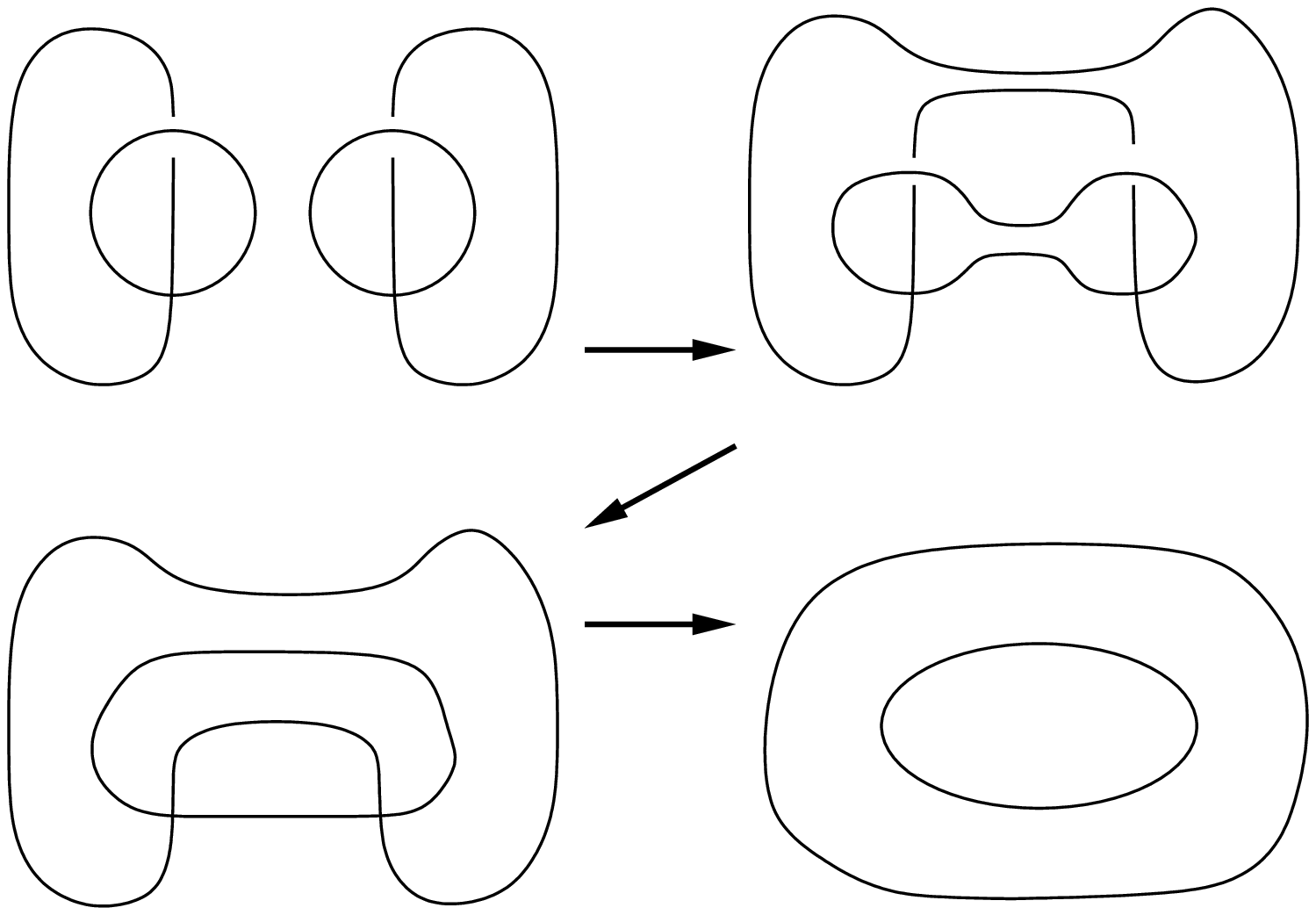}}
\end{center}
\caption{Calceling a pair of pseudo-Hopf links}
\label{hopfcancel}
\end{figure}}

\def\fgKauffvirtualboth{
\begin{figure}[htb]
\begin{center}
\mbox{\epsfxsize=3.0in \epsfbox{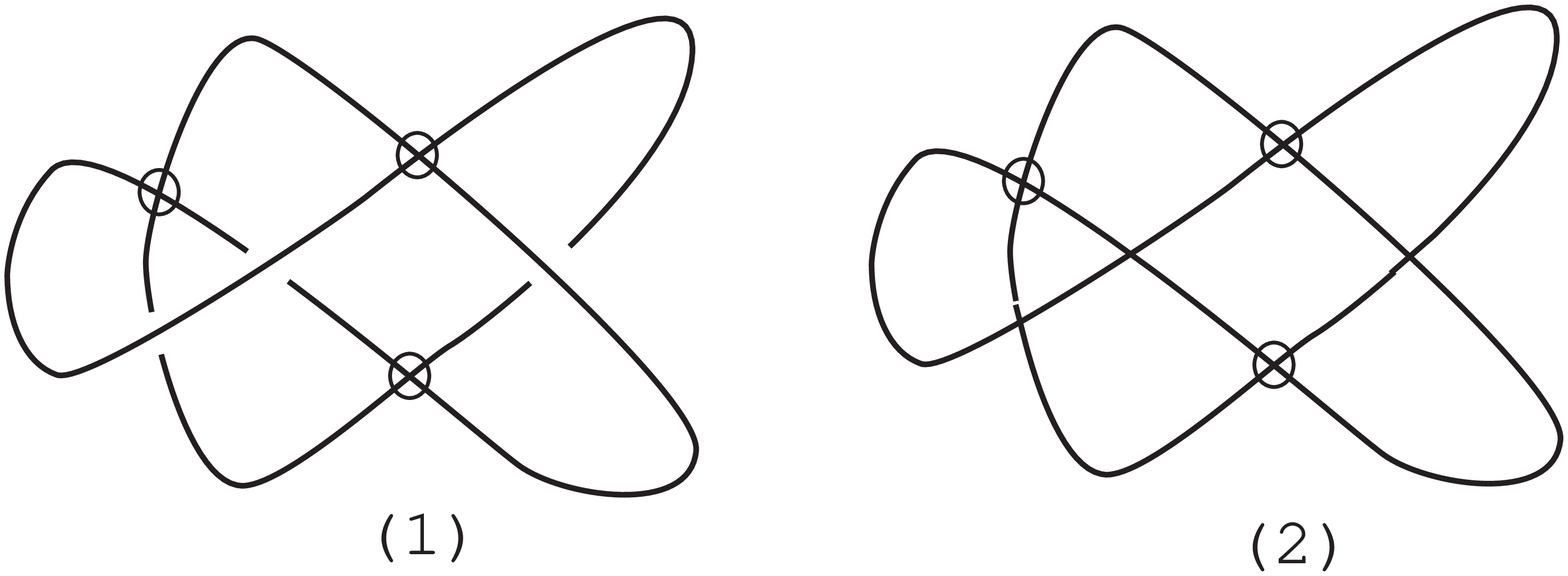}}
\end{center}
\caption{}
\label{Kauffvirtualboth}
\end{figure}}

\def\fglkhom{
\begin{figure}[htb]
\begin{center}
\mbox{\epsfxsize=5in \epsfbox{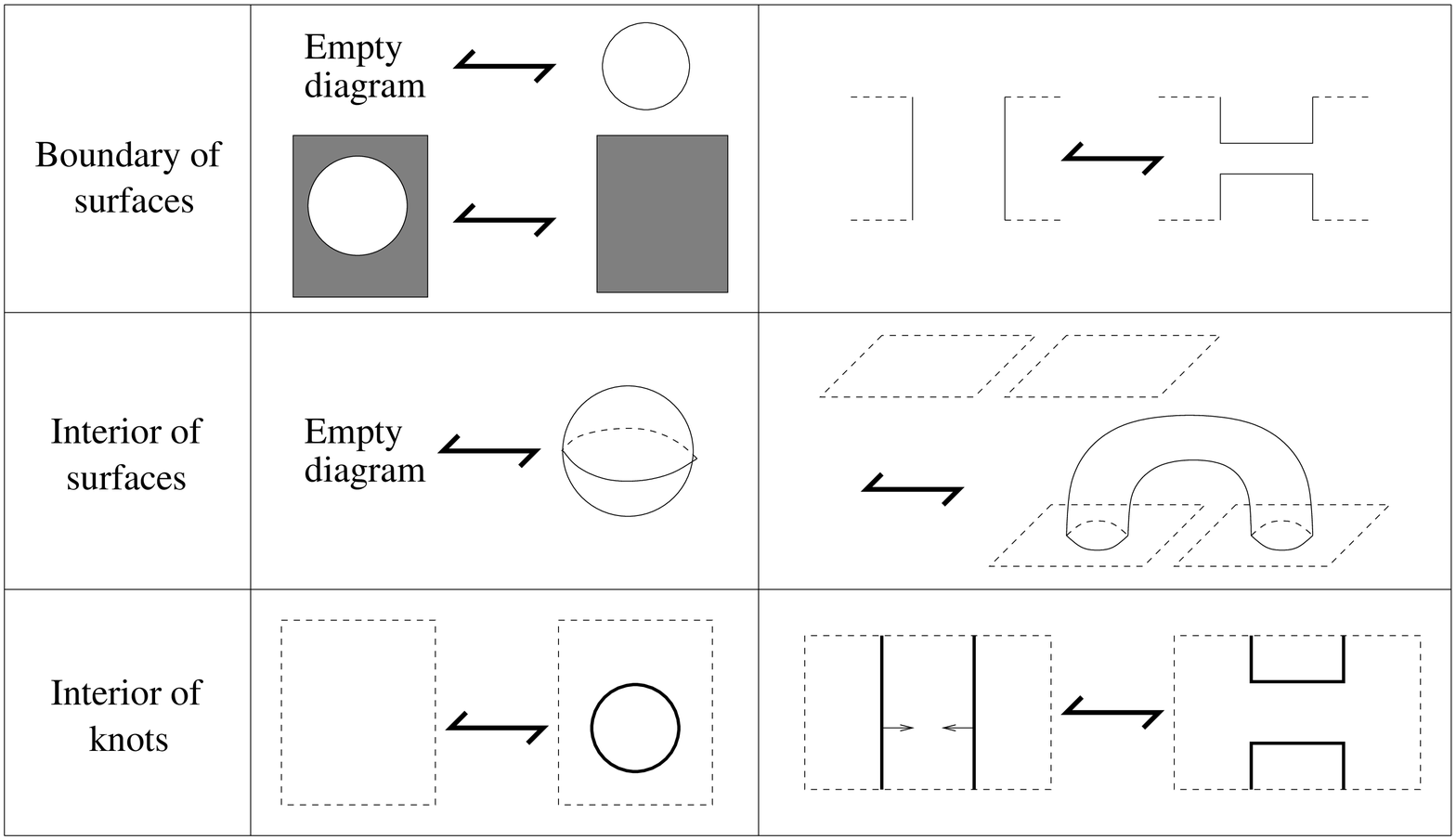}}
\end{center}
\caption{Moves for link-homology}
\label{lkhom}
\end{figure}}

\def\fgskim{
\begin{figure}[htb]
\begin{center}
\mbox{\epsfxsize=5.5in \epsfbox{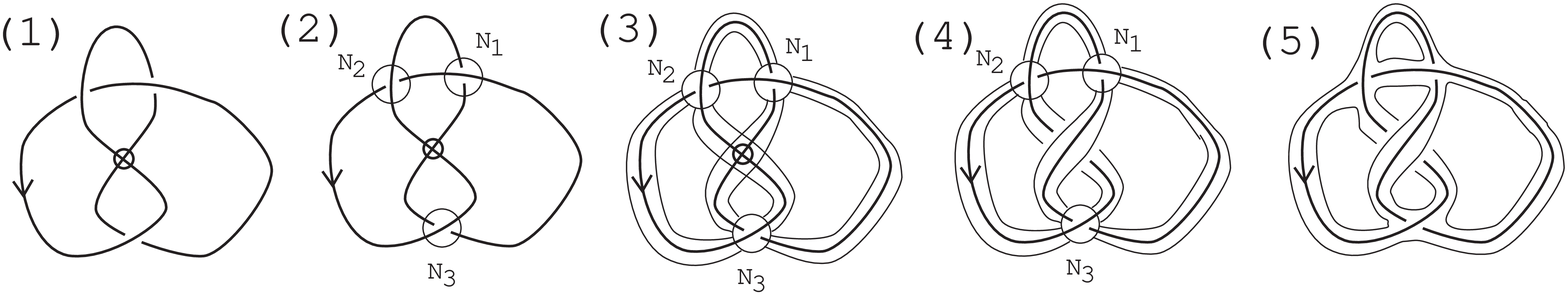}}
\end{center}
\caption{}\label{skim}
\end{figure}}

\def\fgvhopfsp{
\begin{figure}[htb]
\begin{center}
\mbox{\epsfxsize=3.0in \epsfbox{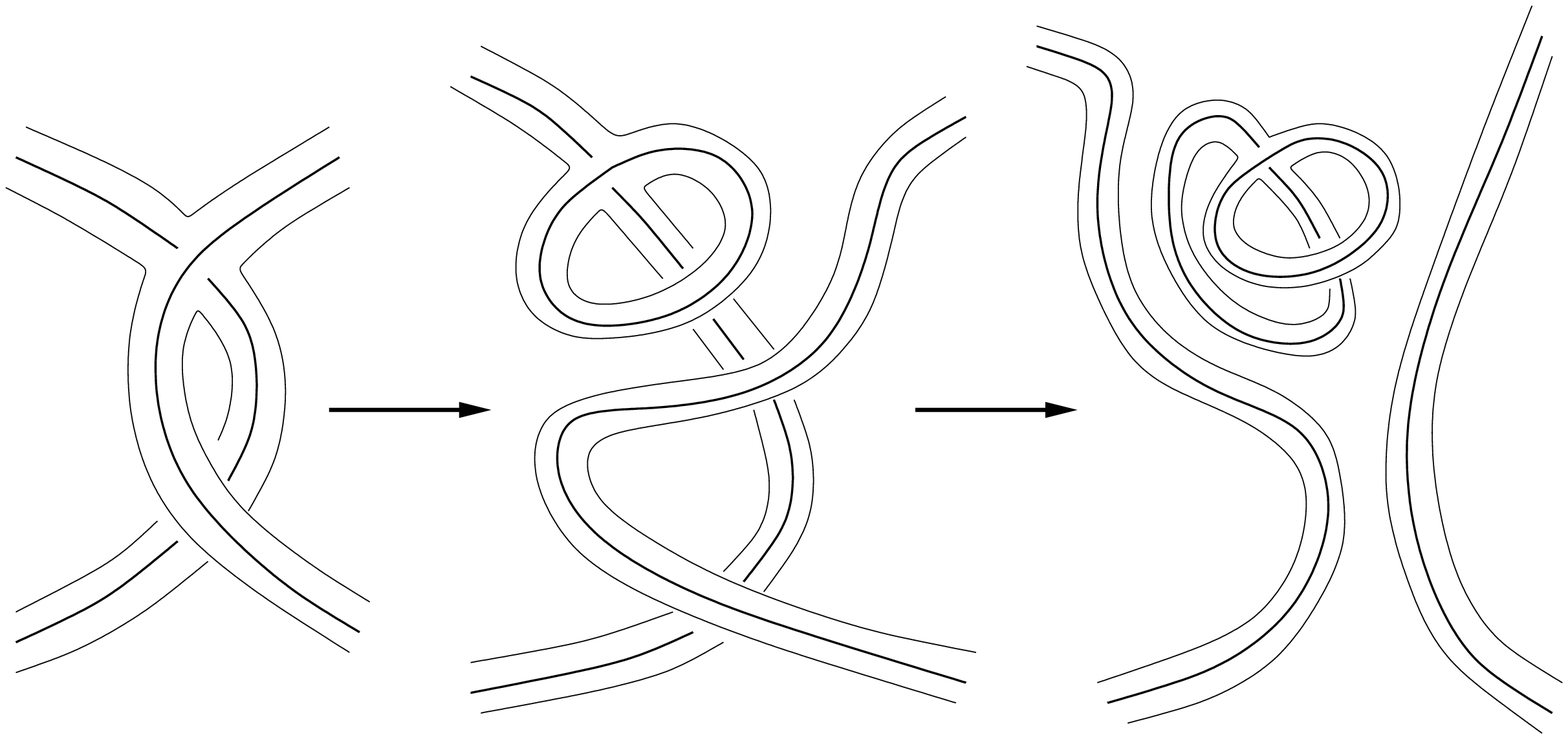}}
\end{center}
\caption{Splitting out a virtual Hopf link}
\label{vhopfsp}
\end{figure}}

\def\fgVRMove{
\begin{figure}[htb]
\begin{center}
\mbox{\epsfxsize=4in \epsfbox{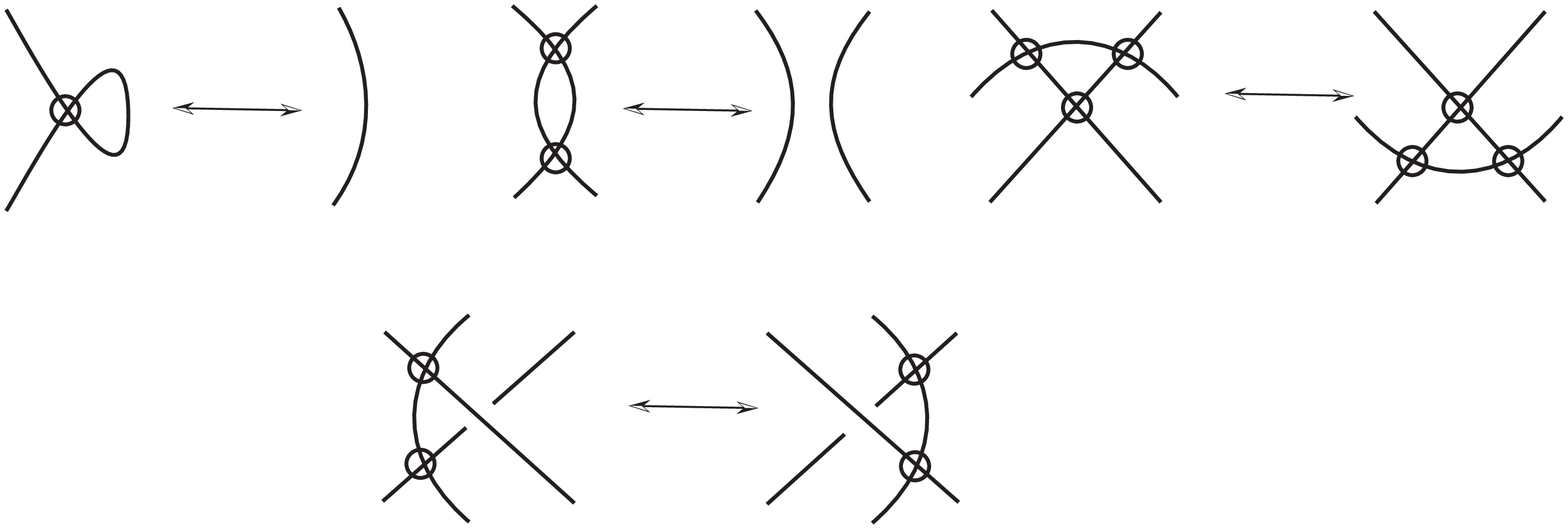}}
\end{center}
\caption{Virtual Reidemeister moves}
\label{VRMove}
\end{figure}}

\def\fgVRMoveA{
\begin{figure}[htb]
\begin{center}
\mbox{\epsfxsize=4in \epsfbox{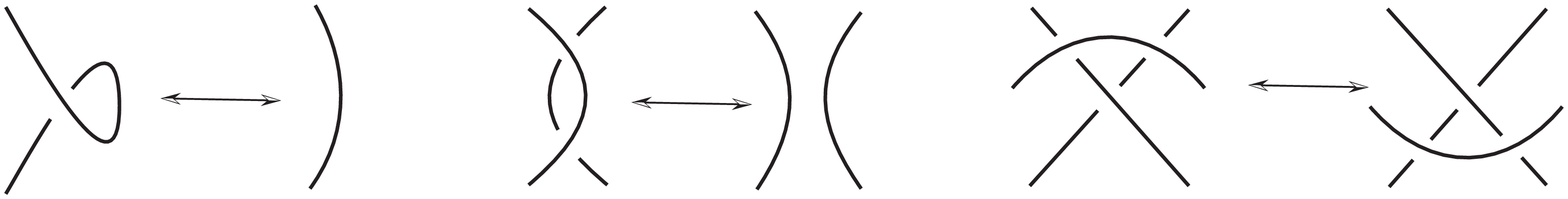}}
\end{center}
\caption{Reidemeister moves}
\label{VRMoveA}
\end{figure}}

\newtheorem{theorem}{Theorem}[section]
\newtheorem{definition}[theorem]{Definition}
\newtheorem{lemma}[theorem]{Lemma}
\newtheorem{corollary}[theorem]{Corollary}
\newtheorem{proposition}[theorem]{Proposition}
\newtheorem{example}[theorem]{Example}
\newtheorem{remark}[theorem]{Remark}
\newcommand{\lk}{\mbox{\rm Link}}
\newcommand{\vlk}{\mbox{\rm vlk}}
\newcommand{\grg}{\mbox{\rm Ground-g}}
\newcommand{\spang}{\mbox{\rm Span-g}}
\newcommand{\suppg}{\mbox{\rm Supp-g}}
\def\qed{$\Box$}

\begin{document}
\title{Stable Equivalence of Knots on Surfaces and
Virtual Knot Cobordisms}

\author{
J. Scott Carter
\footnote{Supported in part by NSF grant DMS-9988107.} \\
University of South Alabama \\
Mobile, AL 36688 \\ carter{\@}mathstat.usouthal.edu
 \and
Seiichi Kamada
\footnote{Supported by 
a Fellowship from the Japan Society
for the Promotion of Science.} \\
Osaka City University \\
Osaka 558-8585, JAPAN\\ kamada{\@}sci.osaka-cu.ac.jp \\
skamada{\@}mathstat.usouthal.edu
\and
Masahico Saito
\footnote{Supported in part by NSF grant DMS-9988101.}\\
University of South Florida \\
Tampa, FL 33620 \\ saito{\@}math.usf.edu
}

\maketitle

\begin{abstract}
We introduce an equivalence relation,
called  stable 
equivalence, on knot diagrams
and closed curves on surfaces.
We give bijections between 
the set of abstract knots, 
the set of virtual knots, and
the set of 
the stable equivalence classes
of knot diagrams on surfaces.
Using these bijections, we define concordance and link homology for
virtual links.
As an application, it is shown that Kauffman's example
of a virtual knot diagram is not equivalent
to a classical knot diagram.

\end{abstract}

\section{Introduction}

Virtual knots were defined in 
\cite{KauB} via diagrams. 
These capture the combinatorial structure of Gauss codes 
and provide interesting 
examples that contrast with classical knot theory.
They were used in \cite{GPV} to study invariants 
of finite type. 
The combinatorial 
nature 
of virtual knots, however,
has caused difficulty in attempts to 
generalize classical invariants.

A bijective 
relation
between virtual knots and certain knots on surfaces,
called abstract knots
was 
given \cite{KK}.
In this paper, we give an alternate geometric interpretation of virtual
knots, called 
stable equivalence of knots on surfaces.
Our interpretation enables us to
introduce notions of cobordisms for virtual knots, for example.
In particular, we classify link homology of virtual links, and use
sliceness to distinguish virtual knots from classical knots
as applications.

The paper is organizes as follows.
In Section~2, we define stable equivalence.
Relations to abstract knots and virtual knots are
established in Section~3.
Cobordisms for virtual knots are defined
and studied in Section~4.
Applications are given in Sections~5.

\section{Stable equivalence 
of knots on surfaces}

Let ${\cal D}$ be the set of all pairs $(F, D)$
such that $F$ is a compact oriented surface and $D$ is a
link diagram on $F$.
   For two elements $(F_1, D_1)$ and $(F_2, D_2)$ of ${\cal D}$,
by
$(F_1, D_1) \overset{e}{\sim} (F_2, D_2)$ we mean that
there exists a compact oriented surface $F_3$
and orientation-preserving
embeddings $f_1: F_1 \to F_3$, $f_2:F_2 \to F_3$
such that
$f_1(D_1)$ and $f_2(D_2)$ are related by Reidemeister moves
on $F_3$ (Fig.~\ref{VRMoveA}).

\fgVRMoveA

\bigskip

\begin{definition} {\rm
{\em Stable Reidemeister equivalence on ${\cal D}$\/} is an
equivalence relation on ${\cal D}$ generated by
the relation $\overset{e}{\sim}$; that is,
two elements $(F, D)$ and $(F', D')$ of ${\cal D}$ are
{\em stably Reidemeister equivalent\/}, denoted by
$(F, D) {\sim} (F', D')$, if there exists a sequence
$(F, D) = (F_1, D_1) \overset{e}{\sim} (F_2, D_2)
\overset{e}{\sim} \cdots \overset{e}{\sim}
(F_n, D_n) = (F', D')$.
}\end{definition}

For example, let $F_1$ and $F_2$ be a torus  $S^1 \times S^1$
and let $D_1$ and $D_2$ be simple closed curves in the torus
such that $D_1$ is null-homotpic and $D_2$ is not.
It is easily seen that $(F_1,D_1) \overset{e}{\sim} (F_2,D_2)$ does not
hold.
  However, $(F_1,D_1) {\sim} (F_2,D_2)$.  Consider an element $(F,D) \in
{\cal D}$ such that $F= S^1 \times [-1,1]$ and $D=S^1 \times \{0\}$.
Then $(F_1,D_2) \overset{e}{\sim} (F,D) \overset{e}{\sim} (F_2,D_2)$.

\bigskip

We note that for an element $(F,D) \in {\cal D}$, a quandle
$Q(D)$ and a group $G(D)$ are defined diagramatically in the usual way in
knot theory.   These are preserved under stable Reidemeister equivalence.

\bigskip

Let ${\cal C}$ be the set of all pairs $(F, C)$
such that $F$ is a compact oriented surface and $C$ is
generic closed curves on $F$. (Generic means that $C$ is
immersed and the 
singularities 
are transverse double points.)
By
$(F_1, C_1) \overset{e}{\sim} (F_2, C_2)$ we mean that
there exists a compact oriented surface $F_3$
and orientation-preserving
embeddings $f_1: F_1 \to F_3$, $f_2:F_2 \to F_3$
such that
$f_1(C_1)$ and $f_2(C_2)$ are homotopic in $F_3$.

\bigskip

\begin{definition} {\rm
{\em Stable equivalence on ${\cal C}$\/} is an
equivalence relation on ${\cal C}$ generated by
the relation $\overset{e}{\sim}$.
}\end{definition}

A natural map
$$ \pi: {\cal D} \to {\cal C}$$
sending a knot diagram to its underlying immersed curve induces a
map
$$ \pi_{\sim}: {\cal D}/_{\sim} \to {\cal C}/_{\sim}.$$
The map $\pi_{\sim}$ is well-defined since homotopy of curves is 
generated by Reidemeister-type 
moves --- more precisely the projection of the Reidemeister moves.

\section{Virtual knots and abstract knots}

\begin{definition} {\rm
(\cite{KauA, KauB}) \quad
A {\em virtual link  diagram\/} consists of 
generic closed curves in
${\bf R}^2$ such that each crossing is either a classical crossing with
over- and under-arcs, or a virtual crossing without over or under
information.
Let ${\cal VL}$ be the set of virtual link diagrams.
The {\em virtual Reidemeister equivalence} is an equivalence relation
on ${\cal VL}$ generated by the Reidemeister moves
depicted  in Fig.~\ref{VRMove}.
Put $VL = {\cal VL}/_{\overset{v}{\sim}}$, where $\overset{v}\sim$ is the
virtual Reidemeister equivalence.  Each element of
$VL$ is called a {\em virtual link\/}.

If the given set of curves of a diagram
is connected (i.e, the diagram consists of
a single component curve), then it is called a {\em virtual
knot 
diagram\/}.
The set of virtual knot diagrams  are denoted by ${\cal VK}$, and
the set of equivalence classes are denoted by
 $VK = {\cal VK}/_{\overset{v}{\sim}}$, whose elements are called
 {\em virtual knots\/}.
}\end{definition}

It is known that there is a bijection between $VK$ and the
set of Gauss codes (or Gauss diagrams) modulo
Reidemeister moves defined in the Gauss code level, \cite{GPV, KauA,
KauB}.  (Refer to \cite{PenneyA, PenneyB} for Gauss codes and
Reidemeister moves on them.)

\fgVRMove

Let ${\cal AL}$ be the subset of ${\cal D}$
consisting of
$(F,D)$ such that $|D|$ is a deformation retract of $F$,
where $|D|$ is the underlying immersed curve in $F$.
See Fig.~\ref{skim}\ (5).  
For $(F_1, D_1), (F_2, D_2) \in {\cal AL}$,
by
$(F_1, D_1) \overset{ae}{\sim} (F_2, D_2)$ we mean that
there exists a {\em closed connected\/} oriented surface $F_3$
and orientation-preserving
embeddings $f_1: F_1 \to F_3$, $f_2:F_2 \to F_3$
such that
$f_1(D_1)$ and $f_2(D_2)$ are related
by Reidemeister moves on $F_3$.

\fgskim

\begin{definition}{\rm
(\cite{NaokoA,NaokoB,NaokoC,KK}) \quad
An {\em abstract link 
diagram} is
an element of ${\cal AL}$.
{\em Abstract Reidemeister equivalence\/},
denoted by $\overset{a}{\sim}$,
is an equivalence relation on ${\cal AL}$
generated by the relation $\overset{ae}{\sim}$.
Put $AL ={\cal AL}/_{\overset{a}{\sim}}$, whose 
elements are 
{\em abstract links \/}. 
}\end{definition}

\begin{theorem} {\rm (\cite{KK})}
There is a map (which we call {\em skimming process})
$$\phi : {\cal VL} \to {\cal AL} $$
that induces a bijection
$$\phi : VL \leftrightarrow AL. $$
\end{theorem}

An abstract link 
diagram is regarded as a
disk-band surface such that there is a usual crossing
in each disk and a proper arc in each band, see
Fig.~\ref{skim}\ (4).
Fig.~\ref{skim} is an illustration of the skimming process,
see \cite{KK} for the definition.

The inclusion map
$$ \iota : {\cal AL} \to {\cal D}$$
induces a map
$$ \iota_{\sim} : AL \to {\cal D}/_{\sim}.$$

\begin{proposition}
The map $ \iota_{\sim} : AL \to {\cal D}/_{\sim}$
is a bijection.
\end{proposition}
{\it  Proof.\/} 
For $(F, D) \in {\cal D}$, let $N(D)$ be a regular neighborhood
of $|D|$ in $F$. Then $(N(D), D)$ is an abstract
link diagram, which we
denote by ${\rm Abs}(F,D)$.
Since ${\rm Abs}(F,D) \overset{e}{\sim} (F, D)$, we see that
the map $ \iota_{\sim}$ is surjective.
Suppose that two abstract
link 
diagrams $(F,D)$ and $(F', D')$
are stably Reidemeister equivalent.
There exists a sequence
$(F, D) = (F_1, D_1) \overset{e}{\sim} (F_2, D_2)
\overset{e}{\sim} \cdots \overset{e}{\sim}
(F_{n-1}, D_{n-1}) \overset{e}{\sim}
(F_n, D_n) = (F', D')$.
Then we have a sequence
$(F, D) = (F_1, D_1) \overset{ae}{\sim} {\rm Abs}(F_2, D_2)
\overset{ae}{\sim} \cdots \overset{ae}{\sim}
{\rm Abs}(F_{n-1}, D_{n-1})  \overset{ae}{\sim}
(F_n, D_n) = (F', D')$.
Thus  $(F, D) \overset{a}{\sim} (F', D')$ and
the map $ \iota_{\sim}$ is injective.  \qed

Now we see that $VL$, $AL$,
${\cal D}/_{\sim}$ and the set of Reidemeister equivalence
classes of Gauss codes are
mutually equivalent.

\section{Link homology and concordance of virtual links}

We recall the definition of knotted surface diagrams \cite{CS:book}.
A knotted surface diagram $K$ is a generically
and properly mapped surface in
a $3$-manifold $M$ such that the double point curves are given
crossing information. Thus $K$ has isolated branch and triple points
and double curves. Along each double curve, one of the two sheets
involved is over-sheet, the other is the under-sheet, and the under-sheet is
broken (interior of small neighborhood removed). At a triple point,
there are top, middle, and bottom sheets. Such a diagram is considered to be
a projection of an embedding of a surface in $M \times [0,1]$.
On the boundary, we have a classical knot diagram on a surface $\partial M$.

\begin{definition} {\rm
Let  $(F_i, D_i)$, $i=0,1$, be two elements of  ${\cal D}$
such that $D_i$ consists of $n$ components $D_i^j$, $j=1, \cdots, n$
for a positive integer $n$.
Then  $(F_0, D_0)$ and $(F_1, D_1)$
are called {\em virtually link-homologous} if there
exists  a compact oriented $3$-manifold $M$ and
a knotted surface diagram $S$
in $M$ with the following properties.

\noindent
(1) $F_0 \cup - F_1 \subset \partial M$, 
where $-F_1$ denotes $F_1$ with its orientation reversed. 

\noindent
(2) $S$ is a knotted surface diagram of an oriented surface
with $n$ components  $S^j$, $j=1, \cdots, n$,
such that $\partial S^j=D_0^j \cup - D_1^j$ for all $j$.
}\end{definition}

\begin{definition} {\rm
Two elements  $(F_i, D_i)$ ($i=0,1$)   of ${\cal D}$ as above,
are called {\em virtually  link-concordant}
(or simply concordant if no confusion occurs) if there exists
a compact oriented $3$-manifold $M$ and a knotted surface diagram $S$
in $M$ with the following properties.

\noindent
(1) $F_0 \cup - F_1 \subset \partial M$.

\noindent
(2) $S$ is a knotted surface diagram of an oriented surface with
$n$ components $S^j$,  $j=1, \cdots, n$, such that
$\partial  C^j=D_0^j \cup - D_1^j$ and each $S^j$ is an annulus.
}\end{definition}

\begin{lemma}\label{linkcobord}
If two elements $(F,D)$ and $(F',D')$ of ${\cal D}$
are stably equivalent, then they are virtually
link-concordant, and hence virtually link-homologous.
\end{lemma}
{\it Proof.\/}
It is sufficient to prove that
if  $(F, D) \stackrel{\mbox{$e$}}{\sim} (F' , D' )$
then   $(F, D)$ and $(F', D')$
are virtually link-concordant.
Let $f: F \rightarrow G$ and $f':  F' \rightarrow G$
be embeddings into a surface $G$
such that $f(D)$ and $f(D')$ are related
by Reidemeister moves in $G$.
Let $M=G \times [0,1]$ and regard $f: F \rightarrow G \times \{ 0 \}$
and $f':  F' \rightarrow G \times \{ 1 \}$, and
identify $F$ and $F'$ with the subsets $f(F)$ and $f'(F')$
of $M$ respectively,
so that $F \cup - F' \subset \partial M$.
Reidemeister moves between $f(F)$ and $f'(F')$ in $G$ yield 
a knotted surface diagram of an annulus in $G \times [0,1]$ such that
the type I, II, and III moves correspond to branch points, 
minimal points of
double point curves, and triple points (cf. \cite{CS:book}),
respectively. 
Hence the result follows.
\qed

\begin{corollary} \label{welldefinedlemma}
The virtual link-concordance and the virtual link-homology
are well-defined for elements of  ${\cal D}/ \sim$.
\end{corollary}

\fglkhom

\begin{lemma} \label{lkhommoveslem}
Let $(F_0, D_0)$ and $(F_1 , D_1)$ be elements of ${\cal D}$.
They are virtually link-homologous if and only if
one is obtained from the other by a sequence of
moves depicted in Fig.~\ref{lkhom}
together with Reidemeister moves.
\end{lemma}
{\it Proof.\/}
Let $(F_0, D_0)$ and $(F_1 , D_1)$ be
virtually link-homologous via a
knotted surface 
diagram
$S$ in a $3$-manifold $M$.
In the following, we regard $S$ as the underlying generic surface
without crossing information for considerations of Morse critical points.
Let $h: M \rightarrow [0,1]$ be a smooth map such that
$h(F_0)=0$ and $h(F_1)=1$.
We may assume (after a small perturbation if necessary)
that $h$
satisfies the following conditions.

(1) $h$ is transverse at $0$ and $1$.

(2) $h$ is generic on $M$, $\partial M$, and $S$,
and on all the self intersections and singularities of $M$, $\partial M$,
and $S$.

Thus $h$ has  isolated Morse
critical points on all the sets listed in (2),
at distinct critical values.
   The sigularities on $S$ gives Reidemeister moves, and
the move listed in Fig.~\ref{lkhom} bottom.
Specifically,
the type I, II, and III moves correspond to
branch points, minimal/maximal 
 points of
double point curves, and triple points.
The minimal/maximal points and saddle points of $S$ 
corresponds to
bottom left and right, respectively,
of Fig.~\ref{lkhom}.
   The Morse critical points as handle moves are
listed in Fig.~\ref{lkhom} top and middle.
The critical points of $\partial M$ are
maxima/minima (the top left entry) 
or saddle points (the top right). 
 {}From the
point of view of the boundary
$1$-manifold, they correspond to handles of
indices $0/2$ and $1$, respectively.
The critical points of Int$M$ are similar,
and depicted in the second row left and right.
Theorem follows as these exhaust generic
singularities and critical points.
\qed

Similarly, we have

\begin{lemma}
Let $(F, D)$ and $(F', D')$ be
elements of ${\cal D}$.
They are virtually link-concordant if and only if
one is obtained from the other by a sequence of Reidemeister moves
and
moves depicted in Fig.~\ref{lkhom},
such that the moves
satisfy the following condition:
the sequence of moves form surface diagrams whose underlying surfaces
are annuli.
\end{lemma}

\begin{definition} {\rm
Let $(F,D)$ be an element of ${\cal D}$ such that
$D$ is a link diagram with $n$ components, $D_j$, $j=1, \cdots, n$.
The {\it linking number\/} between the $j$ and $k$th components,
denoted by $\lk(D_j, D_k)$ is
the number of crossings between $D_j$ and $D_k$ where $D_j$ is
over and $D_k$ is under-arc respectively, counted with signs.
}\end{definition}

This linking number is the same as the
{\it virtual} linking number $\vlk(D_j, D_k)$
in the sense of \cite{GPV}
under the correspondence between virtual links and links on surfaces
via skimming process.

\fghopfboth

\begin{example} {\rm
In the first row of Fig.~\ref{hopfboth},
{\it virtual pseudo-Hopf links\/} are depicted.
The images of them by the skimming process
(in the second row)
are {\it abstract pseudo-Hopf links\/}.
They are {\it positive\/} if the crossings are positive
(right figure);
otherwise {\it negative\/} (left).
The component containing the upper crossing is
called an {\it upper component\/} and the other
a {\it lower component\/}.
Let $D= D_1 \cup D_2$ be a positive
abstract pseudo-Hopf link such that
$D_1$ is upper and $D_2$ is lower, then
$\lk(D_1, D_2) =1$ and $\lk(D_2, D_1) =0$.
}\end{example}

\fgvhopfsp

\begin{proposition}
Virtual link-homology classes of the elements of ${\cal D}$
are completely classified by pairwise linking numbers.
\end{proposition}
{\it Proof.\/}
Bt Lemma~\ref{lkhommoveslem} and the definition,
the linking numbers are invariants of link homology.
We prove the converse.
Since $(F,D) \in {\cal D}$ and ${\rm Abs}(F,D)$
are virtually link-homologous and
have the same linking numbers, we may
assume that $(F,D)$ is an abstract link diagram.
Eliminate each crossing point of $D$
as in Fig.~\ref{vhopfsp}, and we have a split sum of
a trivial abstract link diagram and some
abstract  pseudo-Hopf links.
We cap off each component
of the trivial abstract link diagram.
The remainder is a union of abstract  pseudo-Hopf links.

\fghopfcancel

For each $j$ and $k$ with $j \neq k$,
collect abstract  pseudo-Hopf links
whose upper components come from $D_j$, the $j$th component of
$D$, and the lower components come from $D_k$.
In this family, a pair of positive and negative
abstract  pseudo-Hopf links are canceled as
in Fig.~\ref{hopfcancel}. (For
simplicity,
the figure is
drawn in terms of virtual link diagrams.
Apply the skimming process to obtain
the moves in terms of the abstract  pseudo-Hopf links.)
So we have $| \lk(D_j, D_k)|$  copies of
abstract  pseudo-Hopf links
whose signs are the same with the sign of $\lk(D_j, D_k)$.

Collect abstract  pseudo-Hopf links
whose upper and lower components come from the
$j$th component of $D$, for $j=1,\dots,n$.
If necessary, applying
Reidemeister
moves of type I, 
we may assume that the number of
positive crossings of $D_j$ and the negative crossings
of $D_j$
were
the same.  Then we can eliminate
the abstract  pseudo-Hopf links in this family
as in Fig.~\ref{hopfcancel}.
This implies the proposition.
\qed

\begin{lemma} \label{spanninglemma}
For any $(F, D) \in {\cal D}$, there is
an oriented $3$-manifold $M$ and an oriented surface diagram $G$ in $M$
such that  $\partial M=F$ and  $\partial G= D$.
\end{lemma}
{\it Proof.\/}
Perform a smoothing at each crossing of $D$ to obtain disjoint
simple closed curves $D'$ on $F$.
A smoothing is realized as a branch point. Specifically,
regard $D$ as lying on $F \times \{ 0 \}$ and $D'$ on
$F \times \{ 1 \}$,
then there is an oriented knotted surface diagram $S$
with a branch point corresponding to each smoothing, such that
$\partial S=D \cup -D'$. We cap off each component
of $D'$ by attaching a $2$-handle.
Then we have a desired $M$ and $G$.
\qed

A classical link diagram $D$ on ${\bf R}^2$ is regarded
as an element of ${\cal D}$ by considering $(E, D)$,
where $E$ is a large 2-disk in ${\bf R}^2$
containing $D$ inside.

\begin{corollary}
Two classical link diagrams are virtually link-homologous
if and only if the classical links
represented
by them are link-homologous in classical sense.
\end{corollary}
{\it Proof.\/}
Link homology classes in classical sense are
classified by linking numbers,
whose definition match that of the linking numbers
for elements of ${\cal D}$.
The above proposition, then, implies this corollary.
\qed

\begin{definition} {\rm
An {\it abstractly spanning surface} of $(F, D) \in {\cal D}$ is
a surface $G$ as in Lemma~\ref{spanninglemma}.
The {\it spanning genus} of $(F, D) \in {\cal D}$,
denoted by $\spang(F,D)$, is the minimal genus of all
abstractly spanning surfaces for $(F, D)$.
}\end{definition}

\begin{lemma}
\begin{sloppypar}
If two elements $(F,D)$ and $(F', D')$ are stably equivalent,
then $\spang(F,D) = \spang(F',D')$.
\end{sloppypar}
\end{lemma}
{\it Proof.\/}
This is a consequence of Lemma~\ref{linkcobord}.
\qed

\begin{remark} {\rm
Other kinds of genera of interest are  defined as follows.
A {\it closed realization\/} of $(F, D) \in {\cal D}$
is an embedding of $F$ to a closed oriented surface $G$.
The {\it supporting genus\/} of $(F, D) \in {\cal D}$,
denoted by $\suppg(F,D)$,
is the minimal genus of such closed oriented surfaces $G$,
cf. \cite{NaokoA, NaokoB, NaokoC, KK}.
The {\it ground genus\/} of $(F, D) \in {\cal D}$,
denoted by $\grg(F, D)$, is
the minimal of $\suppg(F',D')$ such that $(F', D')$
is stably equivalent to $(F,D)$.
}\end{remark}

\section{Slice curves on surfaces and Kauffman's example}

Virtual link-concordance on ${\cal D}$ is naturally defined similarly for
${\cal C}$ simply ignoring the crossing informations.
Lemma~\ref{linkcobord} 
holds for ${\cal C}$ under such
a definition, and thus the  virtual link-concordance
is well-defined for ${\cal C}/\sim$.

\begin{definition}
If $(F, D) \in {\cal D}$ or $(F, C) \in {\cal C}$
is virtually concordant to the unlink in the plane,
then it is called {\em slice}.
\end{definition}

By the remark before the definition, we have

\begin{proposition}
Sliceness for ${\cal C}$ is an invariant
under stable equivalence:
Suppose that $(F,C) \sim (F', C')$.
Then $(F,C)$ is slice if and only if $(F', C')$ is slice.
\end{proposition}

In \cite{Carter}, 
a necessary condition for sliceness of immersed closed curves in a surface was given.

\begin{theorem} (\cite{Carter}) \label{Cartersp} 
The pair $(F,C)$ in Fig.~\ref{CarterCurve} is not slice.
\end{theorem}

Figure~\ref{CarterCurve} is different from the example given in
\cite{Carter}.   However it has the same Gauss code with that in
\cite{Carter}  and hence it is not slice.

\fgCarterCurve

In \cite{KauC},
L. Kauffman gave two problems:
\begin{enumerate}
\item
Is the virtual knot diagram in Fig.~\ref{Kauffvirtualboth}\ (1)
virtually Reidemeister equivalent to a classical knot diagram?
(The quandle and the group are the same as those of
a trivial knot diagram.)
\item
Is the universe (Fig.~\ref{Kauffvirtualboth}\ (2)) of the virtual knot
irreducible?
\end{enumerate}
Here a {\em universe\/} of a virtual knot diagram is
a virtual knot diagram without imformation of
over/under crossings for real crossings
(do not confuse them with virtual crossings).
Virtual Reidemeister moves for the universes of virtual knot diagrams
are defined by ignoring over/under information for real crossings.
The universe of a virtual knot is {\it reducible\/}
if it is transformed into the universe of a classical knot
diagram by virtual Reidemeister moves.

\fgKauffvirtualboth

\begin{proposition}
The virtual knot in Fig.~\ref{Kauffvirtualboth}\ (1)
 is not virtually Reidemeister
equivalent to a classical knot diagram.
\end{proposition}
{\it  Proof.\/} 
We have a map
$$\rho =\pi \circ \iota \circ \phi :
{\cal VL} \to
{\cal AL} \to
{\cal D} \to
{\cal C} $$
which induces a map
$$\rho =\pi_{\sim} \circ \iota_{\sim} \circ \phi :
VL  \to
AL  \to
{\cal D}/_{\sim} \to
{\cal C}/_{\sim}. $$
The virtual knot diagram in Fig.~\ref{Kauffvirtualboth}\ (1)
is mapped to $(F,C) \in {\cal C}$ in
Fig.~\ref{CarterCurve}.
This is not slice
(by Theorem~\ref{Cartersp}).  On
the other hand,  any classical knot diagram,
which is regared as an element of ${\cal D}$
by considering it is on a large 2-disk in ${\bf R}^2$, is mapped to
an element of ${\cal C}$ which is slice.
 Since sliceness is 
invariant
under stably equivalence on ${\cal C}$, we see that the virtual knot is
not virtually Reidemeister equivalent to a classical knot
diagram. \qed

Alternate
proofs are given in \cite{Saw} and \cite{SWb}.

\begin{proposition}
The universe in Fig.~\ref{Kauffvirtualboth}\ (2) is not equivalent to
the universe of a classical knot diagram.
\end{proposition}
{\it  Proof.\/} 
The map
$$\pi: {\cal VL} \to {\cal SVL}$$
sending a virtual
link 
diagram to its universe
induces a map
$$\pi_{\sim}: VL \to SVL,$$
where ${\cal SVL}$ is the set of universes of virtual
link 
diagrams
and $SVL$ is the set of equivalence classes.
It is not difficult to see that the map $\rho: VL \to {\cal
C}/_{\sim}$ factors through $SVL$; namely,
when we put $f= \iota_{\sim} \circ \phi$, there is a
map $f'$ which makes the following diagram commutative.
\def\spmapright#1{\smash{
\mathop{\hbox to 1.3cm{\rightarrowfill}}
\limits^{#1}}}
\def\sbmapright#1{\smash{
\mathop{\hbox to 1.3cm{\rightarrowfill}}
\limits_{#1}}}
\def\lmapdown#1{\Big\downarrow
\llap{$\vcenter{\hbox{$\scriptstyle#1\,$}}$}}
\def\rmapdown#1{\Big\downarrow
\rlap{$\vcenter{\hbox{$\scriptstyle#1\,$}}$}}
\begin{equation}
\begin{array}{ccc}
VL & \spmapright{f} & {\cal D}/_{\sim} \\
\lmapdown{\pi_{\sim}~ } & & \rmapdown{\pi_{\sim}} \\
SVL & \sbmapright{f'} & {\cal C}/_{\sim}.
\end{array}
\end{equation}
The universe in Fig.~\ref{Kauffvirtualboth}\ (2)
 is not equivalent to
the universe
of a classical
link 
diagram, because their images
under $f'$ are distinguished
in ${\cal C}/_{\sim}$ by sliceness.
\qed

\end{document}